\documentstyle[german,12pt]{report}
\begin{document}
\noindent
EINE VOLLST\"ANDIGE FORMALISIERUNG DER ARISTOTELISCHEN
NOTWENDIGKEITSSYLLOGISTIK
\par\bigskip
\bigskip
\bigskip
\noindent
I. EINLEITUNG
\par\bigskip
\noindent
Anders als bei den assertorischen Schl\"ussen, die mit einfachen
pr\"adikaten\-lo\-gischen Mitteln nachvollzogen
werden k\"onnen, fehlt es an einer systema\-tischen Interpretation f\"ur die
modale Syllogistik, die Aris\-to\-teles in der ersten Analytik entwickelt.
Dort verwendet er analog zu den vier assertorischen Be\-ziehungen vier
apodiktische Beziehungen, die zwischen zwei Begriffen $A$ und $B$ bestehen k\"onnen, n\"amlich die notwendig allgemein bejahende,
$N(BaA)$, die notwendig allgemein vernei\-nende, $N(BeA)$, die notwendig
partikul\"ar bejahende, $N(BiA)$, sowie die not\-wendig partikul\"ar
vernein\-ende, $N(BoA)$, und untersucht, wann sich aus zwei notwendigen
Pr\"amissen (Kap 8) und wann sich aus einer notwendigen und einer
asserto\-rischen Pr\"amisse (Kap. 9-11) eine notwendige Konklusion ergibt.
\par
Ein formallogischer Zugang zu dieser Theorie der Notwendigkeit hat
die Aufgabe, diese vier Beziehungen in eine strenge Notation zu
\"ubersetzen sowie eine modallogische Axiomatik anzugeben, die eine
m\"oglichst weitgehende \"Ubereinstimmung mit den von Aristoteles
beschriebenen Eigenschaften zu erzielen gestattet. Eine solche Formalisierung
ist im syntaktischen Sinne vollst\"andig zu nennen, wenn sie
simultan folgende f\"unf Anforderungen erf\"ullt.
\par\bigskip
\noindent
(1) Aus den Notwendigkeitsbeziehungen m\"ussen sich die entsprechenden
as\-sertorischen Beziehungen folgern lassen; $N(BaA)$ muss also $BaA$
implizieren (wenn $A$ dem $B$ notwendig zukommt, dann kommt es ihm auch
einfach zu), $N(BeA)$ muss $BeA$ implizieren etc.
\par\bigskip
\noindent
(2) Die notwendig allgemein verneinende und die notwendig partikul\"ar
beja\-hende Beziehung m\"ussen symmetrisch sein, d.h.
$N(BeA) \Longleftrightarrow N(AeB) $ und
$N(BiA) \Longleftrightarrow N(AiB) $.\footnote{Zu dieser Anforderung
vergleiche Aristoteles: Erste Analytik, 25a28ff }
\par\bigskip
\noindent
(3) Aus den notwendig allgemeinen Beziehungen m\"ussen, unter der
zus\"atz\-lichen Voraussetzung, dass ein Begriff nicht leer ist, die
entsprechenden not\-wen\-dig partikul\"aren Beziehungen ableitbar sein, also
$N(BaA) \Longrightarrow N(BiA)$ und
$N(BeA) \Longrightarrow N(BoA)$.\footnote{Aristoteles geht schon
bei den assertorischen Schl\"ussen
Darapti und Felapton davon aus, dass Begriffe niemals leer sind;
zur ersten hier geforderten Implikation s. 25a, die zweite ergibt sich aus
Ferio NXN: $N(BeA) \wedge BiB \Longrightarrow N(BoA) $.}
\par\bigskip
\noindent
(4) Die vierzehn Schl\"usse des achten Kapitels, bei denen zwei notwendige
Pr\"amissen zu einer notwendigen Konklusion f\"uhren, m\"ussen
nachvollziehbar sein.
\par\bigskip
\noindent
(5) Bei den 28 Situationen in den Kapiteln 9-11, wenn der eine
Vordersatz apodiktisch und der andere assertorisch ist, m\"ussen
genau in den von Aristo\-teles beschiebenen dreizehn F\"allen die
Schl\"usse auf Notwendigkeit gehen. F\"ur die \"ubrigen F\"alle m\"ussen
Gegenbeispiele angegeben werden, die zeigen, dass sich hier kein
apodiktischer Schlusssatz ergeben kann.
\par\bigskip
Eine Formalisierung, die im beschriebenen Sinn vollst\"andig ist,
konnte bisher nicht angegeben werden. Die gr\"o\ss ten Schwierigkeiten
ergeben sich hierbei aus dem zuletzt genannten Punkt; auf den ersten
Blick scheint hier Aristoteles je nach Belieben mal einen Beweis 
f\"ur die Notwendigkeit des Schlusssatzes, mal ein Gegenbeispiel
anzuf\"uhren. Schon das erste Paar \--- zwar Barabara NXN,
aber Barbara XNX \--- bereitet Probleme,
und sp\"atestens bei den beiden Aussenseitern
Baroco und Bocardo, bei denen sich, gleichg\"ul\-tig ob die erste oder
die zweite Pr\"amisse notwendig ist, nie ein notwendi\-ger Schlusssatz
ergibt, k\"on\-nen die bis\-heri\-gen For\-mali\-sie\-rungs\-versuche keine
\"Uber\-einstimmung mit der Textvorlage erzielen.
\par
Bemerkenswert ist hierbei, wie die Autoren versuchen, die Abweichungen
ihrer Formalisierungsbem\"uhungen von den aristotlischen Resultaten zu
erkl\"a\-ren. Die einen, die Aristoteles ``wohlgesonnen'' sind,
neigen dazu zu meinen, dass Aristoteles zwar eine klare Vorstellung
von Notwendigkeit habe, dass ihm aber in seinem ansonsten
konsistenten System hin und wieder Fehler unterlaufen seien und er
seinen Ansatz nicht immer konsequent durchgef\"uhrt habe. Die anderen
hingegen nehmen ihre abweichenden Ergebnisse als Beleg f\"ur die
Widerspr\"uchlichkeit der aristotelischen Modallogik und versuchen
darzulegen, dass sie \"uberhaupt nicht konsistent formalisiert werden
und Aris\-toteles sich mit moderner Strenge eben nicht messen lassen
k\"onne.
\par
In dieser Arbeit soll eine vollst\"andige Formalisierung der
Notwendigkeits\-syllogistik des Aristoteles vorgestellt werden. Ihr Ziel
ist damit aber nicht, Aristoteles philosophisch recht zu geben,
sondern lediglich seine modallogi\-schen \"Uberle\-gungen unter Bewahrung
ihrer oft paradox scheinenden Ergeb\-nisse in eine formallogische
Sprache zu \"ubersetzen. Durch Angabe einer kon\-sis\-tenten Interpretation
erbringt sie modelltheoretisch gesprochen den Nachweis f\"ur
die Erf\"ullbarkeit der aristotelischen Modallogik. Ihr Wert be\-steht
somit insbesondere darin, zu zeigen, dass es weder zwingende Gr\"unde
daf\"ur gibt, Aristoteles Widerspr\"uche vorzuwerfen, noch daf\"ur, ihm
Fehler zu unter\-stellen; wenn man m\"ochte kann man sie als Beweis f\"ur
die Aussage lesen, dass die moderne Logik es nicht vermag, die
apodiktische Syllogistik des Aristoteles zu widerlegen.
\par\bigskip
\bigskip
\bigskip
\noindent
II. DIE VIER NOTWENDIGKEITSBEZIEHUNGEN
\par
\bigskip
\noindent
Die Formalisierungen
\par\bigskip
\noindent
Die vier assertorischen Beziehungen, die Aristoteles in den Kapiteln
4-6 der ersten Analytik behandelt, seien pr\"adikatenlogisch formalisiert,
wobei wir gelegentlich auch mengentheoretische Schreibweisen verwenden.
\par\bigskip
\noindent
``$A$ kommt jedem $B$ zu'' durch
\par\medskip
\noindent
$ BaA :\Longleftrightarrow B \subseteq A \Longleftrightarrow
\forall x (Bx \Longrightarrow Ax) \,  ,  $
\par\bigskip\noindent
``$A$ kommt keinem $B$ zu'' durch
\par\medskip
\noindent
$ BeA : \Longleftrightarrow B \cap A = \emptyset \Longleftrightarrow
\forall x ( Bx \Longrightarrow \neg Ax ) \,  ,$
\par\bigskip
\noindent
``$A$ kommt einem $B$ zu'' durch
\par\medskip
\noindent
$BiA : \Longleftrightarrow 
B \cap A  \not = \emptyset \Longleftrightarrow
\exists x (Bx  \wedge Ax ) \,  ,$
\par\bigskip
\noindent
``$A$ kommt einem $B$ nicht zu'' durch
\par\medskip
\noindent
$ BoA : \Longleftrightarrow \exists x( Bx \wedge \neg Ax) \,  .$ 
\par\bigskip
\noindent
Damit lassen sich die assertorischen Schl\"usse problemlos darstellen,
wobei allerdings in der dritten Figur bei Darapti und Felapton darauf
zu achten ist, dass Aristoteles stillschweigend annimmt, dass
Begriffe niemals leer sind.
\par
Diesen vier assertorischen Beziehungen entsprechen, wie schon in
der Einleitung erw\"ahnt, vier apodiktische Beziehungen. F\"ur diese
werden nun modallogische Ausdr\"ucke
angegeben, die die Eigenschaften (1) \--- (5) einer voll\-st\"andigen
Formalisierung erf\"ullen.
\par\bigskip
\noindent
1) Die notwendig allgemein bejahende Beziehung (``$A$ kommt dem $B$
notwen\-dig zu'') sei formalisiert durch
\par\bigskip
\noindent
$N(BaA) : \Longleftrightarrow \forall x( Bx \Longrightarrow
N(Ax)) \,  . $
\par\bigskip
\noindent
2) Die notwendig allgemein verneinende Beziehung (``$A$ kommt
notwendi\-gerweise keinem $B$ zu'', ``$A$ kann kein $B$ sein'')
sei formalisiert durch
\par\bigskip
\noindent
$N(BeA) :\Longleftrightarrow \forall x [Bx \Longrightarrow
\forall y(Ay \Longrightarrow N(x \neq y))] \,  . $
\par\bigskip
\noindent
3)Die notwendig partikul\"ar bejahende Beziehung (``$A$ kommt
notwendig ei\-nem $B$ zu'') sei formalisiert durch
\par\bigskip\noindent
$N(BiA) : \Longleftrightarrow \exists x(Bx \wedge N(Ax))
\, \vee \,  \exists x (N(Bx) \wedge Ax) \, .$
\par\bigskip\noindent
4) Die notwendig partikul\"ar verneinende Beziehung (``$A$ kommt
einem $B$ not\-wen\-dig nicht zu'') sei formalisiert durch
\par\bigskip\noindent
$ N(BoA) :\Longleftrightarrow 
\exists x[Bx \wedge \forall y(Ay \Longrightarrow  N(x \neq y))] $
\par\smallskip
$ \, \, \,  \vee \,   \exists x [N(Bx) \wedge \neg Ax \wedge
\forall y (N(Ay) \Longrightarrow N(x \neq y))] . $\footnote{Die
Formalisierung f\"ur $N(BaA)$ findet sich in einer \"aquivalenten
Form schon bei A. Becker: Die Aristotelische Theorie der
M\"oglichkeitsschl\"usse, M\"unster 1932. Eine zu der von $N(BiA)$ analoge
Formulierung verwendete neulich P. Thom: The two Barbaras, in History
and Philosophy of logic 12 (1991), p.135-149, worauf mich
K.J. Schmidt hingewiesen hat. Die beiden \"ubrigen Formalisierungen
scheinen mir in der Literatur nicht vorzukommen, wie \"uberhaupt
der Gebrauch der notwendigen Verschiedenheit $N(x \neq y)$ neu sein
d\"urfte.}
\par
\bigskip
\bigskip
\noindent
Interpretation in einem modallogischen Modell
\par\bigskip
\noindent
Bevor wir die einzelnen
Formalisierungen kommentieren, wollen wir ein mo\-dal\-logisches Modell angeben,
das es erlaubt, die obigen Formulierungen zu interpretieren. Vom
Grundtyp her wird es sich um ein Leibniz-Kripke-Modell handeln, d.h. die
notwendige G\"ultigkeit wird durch die G\"ultigkeit in einem
Bereich von m\"oglichen Welten erkl\"art. Die auftretenden Terme
$N(Ax)$ und $N(x \neq y)$ erfordern aber zus\"atzliche 
Strukturen in den Modellen, um inter\-pretiert werden zu k\"onnen.
\par\bigskip
\noindent
In dieser Arbeit wird ein modallogisches Modell durch folgende
Daten gege\-ben:
\par\medskip
\noindent
(a) Eine Parametermenge $T$, bei der ein bestimmtes Element $t_0 \in T$
ausge\-zeichnet ist. $T$ kann man als Zeit oder als
M\"oglichkeits\-parameter verstehen, und $t_0$ als Jetztzeit oder als
Wirklichkeitsindex. F\"ur unsere Zwecke gen\"ugt es, $T$ als zweielementig
anzunehmen.
\par\medskip
\noindent
(b) Zu jedem $t \in T$ eine Menge $W_t$, unter der man die Welt
zum Zeitpunkt $t$ verstehen kann. $W_0:=W_{t_0}$ hat man dabei als
Jetztwelt oder Realwelt anzusehen.
\par\medskip
\noindent
(c) Eine Familie von Begriffen $A,B,C $ etc.
Ein Begriff $A$ ist dabei einfach 
eine Familie von Teilmengen $A_t \subseteq W_t,\, t \in T$.
\par\medskip
\noindent
(d) Eine Familie von Individuen oder besser Individualkonzepten $x,y $ etc.
Ein Individualkonzept $x$ ist dabei eine Abbildung
\par\smallskip
\noindent
$x :T \longrightarrow W:= \biguplus_{t \in T} W_t $ mit $x_t=x(t) \in W_t$.
\par\smallskip
\noindent
Ein Individualkonzept ordnet also jedem Parameter $t$ ein Element in der
zu $t$ geh\"orenden Welt $W_t$ zu. Man verwechsele ein Individuum $x$ nicht
mit seiner realen Erscheinungsweise $x_0 \in W_0$!
\par\bigskip
\noindent
Ist ein solches Modell gegeben, so lassen sich die Ausdr\"ucke in
nahe liegender Weise interpretieren. Nicht durch Quantoren gebundene Terme
setzen f\"ur ihre Interpretation eine Belegung der Individualvariablen und
der Begriffsva\-riablen voraus.
\par\medskip
\noindent
$Ax$  steht f\"ur  $x_0 \in A_0 $.
\par
\medskip
\noindent
$N(Ax)$  steht f\"ur  $ x_t \in A_t $ f\"ur alle $t \in T$.
\par\medskip
\noindent
$N(x \neq y)$  steht f\"ur  $x_t \neq y_t $ f\"ur alle $t \in T$.
\par\bigskip
\noindent
Quantifiziert wird stets \"uber die Individuen, nicht etwa \"uber die
Realwelt $W_0$. Die verwendeten Formalisierungen lassen sich in dieser
Weise in einem solchen Modell interpretieren und auf ihre
G\"ultigkeit \"uberpr\"ufen. Auch verschiedene andere, in der Literatur
zum Verst\"andnis der aristotelischen Modallogik vorgeschlagenen
Ausdr\"ucke lassen sich derartig interpretieren und mit den hier
vertretenen Formeln vergleichen. Unabh\"angig von der antiken Vorlage
erlaubt das Modell, verschiedene Arten des notwendigen Zukommens
zu durchdenken und deren Verh\"altnis zu kl\"aren.\footnote{Analytische
Beziehungen wird man wohl so anzusetzen haben, 
dass die zu Grunde liegende assertorische Beziehung in allen Welten allein
unter Bezug auf den Begriff, und nicht unter Bezug auf das m\"ogliche
Verhalten der
in einer Startwelt anzutreffenden Individuen. Die analytische Inklusion
w"are demnach mit $\forall t (A_t \subseteq B_t) $ zu fassen; das 
Junggesellen unverheiratet sind, ist von diesem Typ. Junggesellen
sind aber nicht notwendig unverheiratet, ebensowenig wie sie notwendig
Junggesellen sind. Das w\"urde bedeuten, dass die Junggesellen von heute
Zeit ihres Lebens Junggesellen bleiben m\"ussten, was allenfalls f\"ur
eingefleischte Junggesellen gilt.}
Die modallo\-gischen Modelle werden in dieser Arbeit in erster Linie dazu
verwendet, Gegenbeispiele zu den gemischtmodalen Situationen zu
konstruieren, bei denen nach Aristoteles kein Schluss auf
Notwendigkeit stattfindet.
\par\bigskip
\noindent
Die modallogischen Formulierungen werden mittels den Modellen
auf pr\"adi\-katenlogische For\-mulierungen zur\"uckgef\"uhrt, und die Nachzeichnung der
ari\-stotelischen Schl\"usse wird sich daher allein pr\"adikatenlogischer
Mittel bedie\-nen, ohne weiterer modallogischer Gesetze im eigentlichen Sinn.
Auch Aris\-toteles f\"uhrt seine apodiktischen Schl\"usse 
letztendlich auf assertorische Be\-zieh\-ungen zur\"uck, zumindest legen das
Redewendungen wie ``denn $C$ steht unter $A$'',
``da $A$ ein $B$ ist'',
``nun kommt aber $A$ einem $C$ zu'' nahe, die letztlich in den Kapiteln 9-11
die Schl\"usse begr\"unden.
\par
\bigskip
\bigskip
\noindent
Erl\"auterungen zu den einzelnen Formalisierungen
\par\bigskip
\noindent
Die vorgestellten Formulierungen f\"ur die vier modallogischen
Be\-zie\-hung\-en 
des Aris\-to\-te\-les bed\"urfen einiger Erl\"auterungen, die
gemein\-sam mit dem Nach\-weis der elementaren Eigenschaften (1) \--- (3)
gegeben werden sollen.
\par\bigskip
\noindent
Die Formalisierung von
$N(BaA)\Longleftrightarrow \forall x(Bx \Longrightarrow N(Ax))$
besagt, dass jedes Individuum, das zur 
Jetztzeit zu $B$ geh\"ort (real unter den Begriff $B$ f\"allt),
mit Notwendigkeit zu $A$ geh\"ort, also nicht nur im realen Sinn, sondern
auch noch in den variierten Welten. Er bedeutet also
$\forall x( x_0 \in B_0 \Longrightarrow \forall t(x_t \in A_t))  .$   
\par
Da die notwendige Zugeh\"origkeit $NAx$ die reale Zugeh\"origkeit
$Ax \Longleftrightarrow x_0 \in A_0 $ impliziert, kann man aus
der G\"ultigkeit von $N(BaA)$ sofort auf die von $BaA$ schlie\ss en.
Unter der zus\"atzlichen Voraussetzung, dass $B$ nicht
leer ist, (d.h. $ B_0 \neq \emptyset $)
folgt aus $N(BaA)$ die Existenz eines $x \in B$, f\"ur das
$NAx$ gilt, d.h. es l\"asst sich der Ausdruck
$\exists x (Bx \wedge NAx)$ gewinnen,
womit die G\"ultigkeit der Implikation $N(BaA) \Longrightarrow N(BiA)$
nachgewiesen ist.
\par\bigskip
\noindent
Der Ausdruck f\"ur $N(BeA) \Longleftrightarrow
\forall x (Bx
\Longrightarrow \forall y[Ay \Longrightarrow N(x \neq y))]$
besagt, dass jedes Element aus $B$ von jedem
Element aus $A$ notwendigerweise ver\-schieden ist, dass also zwei Individuen,
wo real das eine zu $A$ und das andere zu $B$ geh\"ort,
niemals gleich sein (oder werden) k\"onnen.
Eine einfache 
pr\"adika\-tenlogische Umformulierung ergibt
$$ N(BeA) \Longleftrightarrow
\forall x \forall y
(Bx \wedge Ay \Longrightarrow N(x \neq y))\,  ,$$
woran man sofort die Symmetrie in $A$ und $B$ ablesen kann.
Die Formalisie\-rung erlaubt also die \"Aquivalenz
$N(BeA) \Longleftrightarrow N(AeB)$.
Ohne Notwendig\-keits\-operator l\"asst sich $N(BeA)$ schreiben als
$$\forall x \forall y \forall t
(x_0 \in B_0 \wedge y_0 \in A_0 \Longrightarrow(x_t \neq y_t)) \,   .$$
Daraus gewinnt man sofort auch die assertorische Disjunktheit von
$A$ und $B$.
\par
Unter der zus\"atzlichen Voraussetzung, dass $B$
nicht leer ist, sagen wir $x \in B$, folgt aus der notwendigen Disjunktheit
$N(BeA)$ die Existenzaussage 
$\exists x [Bx \wedge \forall y(Ay \Longrightarrow N(x \neq y)]$
und damit $N(BoA)$.
\par\bigskip
\noindent
Aufgrund der Symmetrie von $\vee$ und $\wedge$ sind $N(BiA)$ und $N(AiB)$
\"aquiva\-lent. Der Ausdruck f\"ur $N(BiA)$ bedeutet,
dass es ein Element $x \in B \cap A$ gibt, das zus\"atzlich in mindestens
einer der Mengen $A$ und $B$ mit Notwendigkeit enthalten ist.
\par\bigskip
\noindent
Die Formalisierung f\"ur $N(BoA)$ scheint auf den ersten Blick ziemlich
kompli\-ziert; das folgende Implikationsdiagramm soll die Situation
verdeut\-lichen.
\par\bigskip
\noindent
\unitlength0.5cm
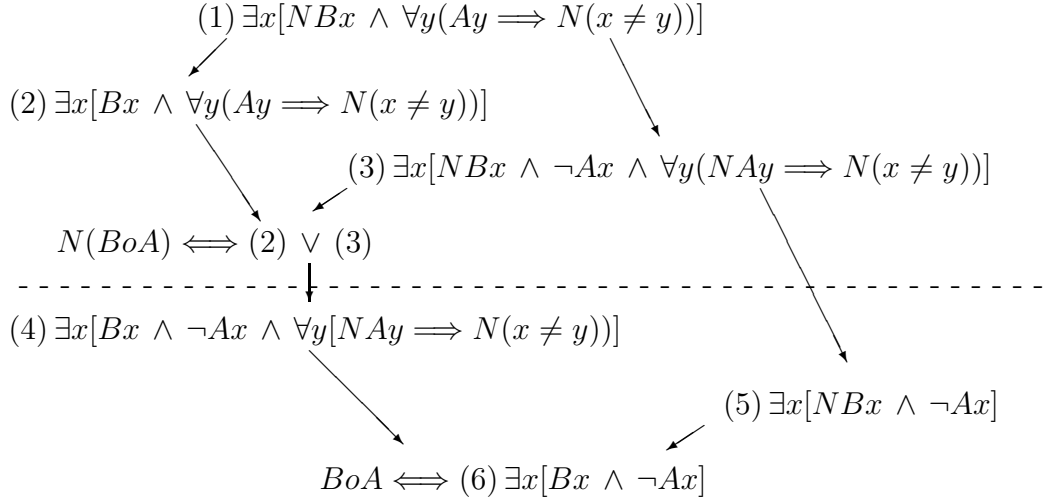
\begin{picture}(25,17)
\put(6,15.8){\rm Implikationsdiagramm f\"ur $N(BoA)$}
\put(5,13.3){$ (1)\,  \exists x[NBx \, \wedge \,
\forall y (Ay \Longrightarrow N(x \neq y))]$}
\put(5.8,13){\vector(-1,-1){1}}
\put(0,11) {$ (2) \, \exists x [Bx \, \wedge \,
\forall y(Ay \Longrightarrow N(x \neq y))]$}
\put(16,13){\vector(1,-2){1.3}}
\put(9,9.3){$ (3)\, \exists x[NBx \, \wedge \, \neg Ax \, 
\wedge \, \forall y ( NAy \Longrightarrow  N(x \neq y))] $}
\put(5,10.7){\vector(2,-3){1.7}}
\put(9,9){\vector(-3,-2){0.9}}
\put(20,9){\vector(1,-2){2.3}}
\put(1,7.3) {  $ N(BoA) \Longleftrightarrow (2) \, \vee \, (3) $}
\put(0,6.2){\- - - - - - - - - - - - - - - - - - - - - - - - - - -
- - - - - - - - - - - - - - - - - - - - - - - - }
\put(8,7){\vector(0,-1){1}}
\put(0,5) {$ (4) \, \exists x[Bx \, \wedge \, \neg Ax \, \wedge \,
\forall y [NAy \Longrightarrow N(x \neq y))] $}
\put(19,3){$ (5) \, \exists x[NBx \, \wedge \, \neg Ax]$}
\put(8,1){ $ BoA \Longleftrightarrow (6)\, \exists x[Bx \,
\wedge \, \neg Ax]$}
\put(8,4.7){\vector(1,-1){2.5}}
\put(18.5,2.7){\vector(-3,-2){1}}
\end{picture}
\par\bigskip
\noindent
Diese Implikationen beruhen allesamt auf dem $T-$Axiom und aussagenlo\-gi\-schen
Gesetzen.\footnote{$NBx \Longrightarrow Bx$ liefert
$(1) \Longrightarrow (2), (3) \Longrightarrow (4)$ und
$(5) \Longrightarrow (6)$.
$ \forall y (Ay \Longrightarrow N(x \neq y))$ impliziert einerseits
$ \neg Ax$, andererseits
$\forall y (NAy \Longrightarrow N(x \neq y))$,
was die Implikationen $(1) \Longrightarrow (3)$ und $(2) \Longrightarrow (4)$
gew\"ahrleistet.}
Die gestrichelte Linie markiert die Grenze
dessen, was im Rahmen dieser Interpretation noch als notwendig
partikul\"ar verneinend angesehen wird und was nicht mehr. Das Diagramm
besagt mittels $(6)$, dass $BoA$ aus $N(BoA)$ folgt. Die Stellung
von (5) l\"asst erkennen, dass es keinen Zusammenhang zwischen $N(BoA)$
und $(N(Bi \neg A)$ gibt. (4) wird beim Beweis von Baroco NNN und Bocardo NNN
eine Rolle spielen. Dass (2) als notwendig partikul\"ar verneinend
zu gelten hat folgt unmittelbar aus der Formalisierung f\"ur
$N(BeA)$ und Ferio NXN. W\"urde man jedoch $N(BoA)$ einfach gleich (2)
ansetzen, so ergebe sich bei Baroco XN und Bocardo NX ein
Notwendigkeitsschluss im Gegensatz zu Aristoteles Behauptungen in den
Kapiteln 10 und 11. L\"asst man aber zus\"atzlich noch (3) als
notwendig partikul\"ar verneinend gelten, so findet dort ein solcher 
Schluss nicht mehr statt. Indem $N(BoA)$ als die alternative
Verkn\"upfung der Ausdr\"ucke (2) und (3) angesetzt wird, erh\"alt man wie
sich zeigen wird eine Formulierung, die einerseits schwach genug
ist, um in den Kapiteln 10 und 11 keinen Schluss zuzulassen, andererseits
aber stark genug, um in Kapitel 8 Baroco NNN und Bocardo NNN
als g\"ultig zu erweisen.
\par
Da das notwendig partikul\"ar verneinende bei den
apodiktisch\--assertori\-schen Situationen nie zu einem notwendigen
Schlusssatz f\"uhrt, k\"onnte man auch auf eine explizite Formalisierung
f\"ur $N(BoA)$ ganz verzichten, und stattdessen einfach
$\exists x[Bx \wedge \forall y(Ay \Longrightarrow N(x \neq y))]
\Longrightarrow N(BoA) $
sowie die Schl\"usse Baroco NNN und Bocardo NNN direkt als
Axiome postulieren. Das scheint mir aber mindestens so kompliziert
und auf jeden Fall unbefrie\-digender zu sein als der vorgestellte
Ansatz. In irgendeiner Weise muss eine formale Interpretation die
Komplexit\"at der Vorlage beinhalten: entweder in einem aufgebl\"ahten
Regelapparat oder in relativ komplizierten Formalisie\-rungen. Ich
verfolge hier die zweite Strategie. So erkl\"aren sich beispielsweise
die Formalisierungen zu $N(BeA)$ und $N(BiA)$ unter anderem dadurch,
die Sym\-metrieeigenschaft bereits in den Formulierungen zu sichern,
ohne hierzu eigene Axiome aufzustellen mit dem Risiko, Schl\"usse zu
erhalten, die nach Aristoteles ung\"ultig sind.
\par\bigskip
\bigskip
\noindent
Bemerkungen zu der Interpretation
\par
\bigskip
\noindent
Zu der vorgestellten formalen Interpretation als Ganzes ist zu
sagen, dass sie gewisse, von einem naiven Verst\"andnis\footnote{Damit
meine ich ein Verst\"andnis, das sich allzu sehr von den
assertorischen Bege\-ben\-heiten leiten l\"asst.}
her zu erwartende
Eigenschaften gerade nicht aufweist. Einige dieser Eigenschaften
wollen wir hier herausar\-beiten und zeigen, dass sie von der aristotelischen
Vorlage her auch nicht zu er\-war\-ten sind.
\par
Die notwendig
verneinenden Beziehungen gehen nicht aus den entsprech\-enden bejahenden
Beziehungen hervor, indem man an geeigneter Stelle eine Ne\-ga\-tion
einf\"uhrt, d.h. es gilt weder $N(BeA) \Longleftrightarrow N(B a \neg A)$ noch
$N(BoA)$ 
\par\noindent
$\Longleftrightarrow N(B i \neg A)$.
Damit in Zusammenhang steht die Ung\"ultigkeit der ``apo\-diktischen
Kontraposition''  $ N(BaA) \Longleftrightarrow N(\neg A a \neg B) $.
Die erste
\"Aquivalenz verbietet sich, wenn man Barbara XNX betrachtet: dessen
Pr\"amissen $BaA$ und $N(CaB)$ w\"aren dann n\"amlich \"aquivalent zu
$N(\neg B e C)$ und $\neg A a \neg B $,
was mittels Celarent NXN die Konklusion $N(\neg A e C)$, also $N(CaA)$
ergeben w\"urde. Die zweite \"Aquivalenz kann zumindest nicht gemeinsam
mit der apodiktischen Kontraposition gelten, ansonsten lie\ss en sich die
Pr\"amissen von Baroco NXX, $N(BaA)$ und $CoA$, zu $N(\neg A a \neg B)$ und
$C i \neg A$ umformulieren, was mit Darii NXN zu $N(Ci \neg B)$, also $N(CoB)$
f\"uhren w\"urde.
\par
Des weiteren entstehen die partikul\"aren Beziehungen nicht aus den
all\-ge\-meinen, indem man einfach einen Allquantor $\forall$ durch einen
Existenz\-quantor $\exists $ ersetzt. Wendet man dieses Verfahren auf die
Formalisierungen f\"ur $N(BaA)$ bzw. $N(BeA)$ an, so ergibt sich jeweils
nur die eine H\"alfte der For\-malisie\-rungen f\"ur $N(BiA)$ bzw. $N(BoA)$.
\par
Daneben ergeben sich die apodiktischen Beziehungen aus den assertori\-schen
nicht dadurch, dass diese als Ganze mit einem Notwendigkeitsoperator $N$
versehen werden. Insbesondere sind Schreibweisen wie $N(BaA)$, $N(BeA)$ etc.
lediglich als Abk\"urzungen f\"ur die modalen Beziehungen 
zu verstehen, nicht aber als $N(\forall x(Bx \Longrightarrow Ax))$,
$N(\forall x(Bx \Longrightarrow \neg Ax))$ etc.
\par
Weiterhin ist bemerkenswert, dass Tautologien im allgemeinen nicht
not\-wendig im aristotelischen Sinne sind. W\"urden n\"amlich
beispielsweise aus den Tautologien $AaA$, d.h. alle $A$ sind $A$, und
$\neg A e A$, d.h. alle Dinge sind $A$ oder nicht $A$, die
Notwendigkeitsbeziehungen $N(AaA)$ und $N(A e \neg A)$  f\"ur beliebiges
$A$ folgen, so ergebe sich aus $BaA$ (bzw BiA bzw. BeA)
sofort mittels Barbara NXN (bzw. Darii NXN bzw. Celarent NXN)
die G\"ultigkeit von N(BaA)
(bzw. die von N(BiA) bzw. N(BeA))\footnote{Barbara NXN:
$N(AaA) \wedge BaA \Longrightarrow N(BaA)$,
Darii NXN: $N(AaA) \wedge BiA \Longrightarrow N(BiA)$,
Celarent NXN: $N(\neg A e A) \wedge B a \neg A \Longrightarrow N(B e A)$
($\neg A $ als Zwischenbegriff). }.
Dies w\"urde nicht nur Schl\"usse wie Barbara XNN, Celarent XNN etc.
gestatten, sondern auch \--- wenn man das $T-$Axiom gelten l\"asst und die
Schl\"usse als Implikationen deutet \--- die \"Aquivalenzen
$N(BaA) \Longleftrightarrow BaA$, $N(BiA) \Longleftrightarrow BiA$ und
$N(BeA) \Longleftrightarrow BeA$ ergeben, was die gesamte
aristotelische Notwendigkeits\-syllogistik unsinnig machen w\"urde.
\par\bigskip
\bigskip
\bigskip
\noindent
III. DIE NOTWENDIGKEITSSCHL\"USSE DES ARISTOTELES
\par
\bigskip
\noindent
Es ist nun zu zeigen, dass die oben vorgestellte Formalisierung auch
die Anforderungen (4) und (5) erf\"ullt, wobei aus beweistechnischen
Gr\"unden (5) vorgezogen wird. Dabei werden je Figur die aristotelischen
Schl\"usse bewiesen und anschlie\ss end f\"ur die Situationen, in denen
nach Aristoteles keine Schl\"usse auf Notwendigkeit stattfinden,
Gegenbeispiele angegeben.
\par
\bigskip
\bigskip
\noindent
1.Figur
\par\bigskip
\noindent
\unitlength0.5cm
\begin{picture}(10,6)
\put(0,5.3){\rm Barbara\ NXN }
\put(0,4){$ \forall x (Bx \Longrightarrow NAx) $}
\put(0,3){$ \forall x(Cx \Longrightarrow Bx) $}
\put(0,2.6){\----------------------------}
\put(0,1.7){$ \forall x(Cx \Longrightarrow NAx)$}
\put(13,5.3){\rm Celarent NXN}
\put(13,4){$\forall x \forall y (Bx \wedge Ay \Longrightarrow N(x \neq y)) $}
\put(13,3){$ \forall x(Cx \Longrightarrow Bx) $}
\put(13,2.6){\---------------------------------------------}
\put(13,1.7){$ \forall x \forall y(Cx \wedge Ay
\Longrightarrow N(x \neq y))) $}
\end{picture}
\par\bigskip
\noindent
\unitlength0.5cm
\begin{picture}(10,6)
\put(0,5.3){\rm Darii\ NXN }
\put(0,4){$ \forall x (Bx \Longrightarrow NAx) $}
\put(0,3){$ \exists x(Cx \wedge Bx) $}
\put(0,2.6){\------------------------------}
\put(0,1.7){$ \exists x(Cx \wedge NAx)$}
\put(0,0.4){\rm also $N(CiA)$. }
\put(13,5.3){\rm Ferio NXN}
\put(13,4){$\forall x[Bx \Longrightarrow
\forall y (Ay \Longrightarrow N(x \neq y))] $}
\put(13,3){$ \exists x(Cx \wedge Bx) $}
\put(13,2.6){\------------------------------------------------}
\put(13,1.7){$ \exists x[Cx \wedge
\forall y(Ay \Longrightarrow N(x \neq y))] $}
\put(13,0.4){\rm also $N(CoA)$. }
\end{picture}
\par\bigskip
\noindent
Dagegen findet in den \"ubrigen Situationen kein Schluss auf
Notwendigkeit statt.\footnote{Findet kein Notwendigkeitsschluss
statt, so wird, wenn die Formalisierung eine Alternation aufweist
(also bei $N(BiA)$ und $N(BoA)$) nur der ung\"unstige Fall aufgef\"uhrt.} 
\par\bigskip   
\bigskip
\noindent
\unitlength0.5cm
\begin{picture}(10,4)
\put(0,3.3){\rm Barbara XN }
\put(0,2){$ \forall x (Bx \Longrightarrow Ax) $}
\put(0,1){$ \forall x(Cx \Longrightarrow NBx) $}
\put(13,3.3){\rm Darii XN}
\put(13,2){$ \forall x(Bx \Longrightarrow Ax) $}
\put(13,1){$ \exists x(Cx \wedge NBx) $}
\end{picture}
\par
\noindent
Wie konstruieren hierzu ein gemeinsames Gegenbeispiel, indem wir 
ein Mo\-dell angeben, indem die Pr\"amissen von Barbara XN (und von
Darii XN) er\-f\"ullt sind, aber nicht die in Frage stehende Konklusion.
Es sei $T= \{ 0 ,1 \} $ mit $W_0 =\{ x_0 \},\, W_1 = \{ x_1 \} $. 
Das einzige Individuum ist nat\"urlich $x=(x_0,x_1)$.
Die Begriffe $A,B,C$ seien gegeben durch
$A_0=B_0=C_0=\{ x_0 \} $ und $A_1=C_1= \emptyset \,  , B_1 = \{ x_1 \} \, .$
Damit ist $NBx$, also auch $C \subseteq NB $. Dagegen ist
$\neg NAx $ und $\neg NCx $ und auch nicht
$C \subseteq NA$.
Als Konklusion ergibt sich nur die Existenz eines
gemeinsamen Elements von $A$ und $C$, ohne dass dieses in einer dieser
Mengen mit Notwendigkeit enthalten ist.
\par
Aristoteles selbst f\"uhrt die Annahme, es ergebe sich bei Barbara XN ein
Schluss auf Notwendigkeit, zum Widerspruch: ``Denn w\"are der
Schlusssatz notwendig, so folgte sowohl durch die erste Figur wie auch
durch die dritte Figur, dass $A$ einem $B$ notwendig zukommt. Das ist aber
falsch. Denn $B$ kann so beschaffen sein, dass m\"oglicherweise $A$
keinem $B$ zukommt.''\footnote{Erste Analytik, 30a25ff}
Auch das l\"asst sich formal
nachzeichnen. Aus der Annahme $\forall x(Cx \Longrightarrow NAx $
folgt n\"amlich (unter der zus\"atzlichen Voraussetzung $C \neq \emptyset$)
zun\"achst $\exists x(Cx \wedge NAx)$ und daraus mit der zweiten
Pr\"amisse $\exists x(NBx \wedge NAx)$.
Die erste Pr\"amisse $\forall x (Bx \Longrightarrow Ax)$ ist aber
im Widerspruch hierzu durchaus mit
$\forall x(Bx \Longrightarrow M(\neg Ax))$\footnote{Gelegentlich
verwende ich den M\"oglichkeitsoperator $M$; es gilt die Beziehung
$ \neg Np \Longleftrightarrow M \neg p$. Insbesondere ist also
$M(\neg Ax) \Longleftrightarrow \neg N(Ax)$ und
$ \neg N(x \neq y) \Longleftrightarrow M(x=y) $.}
vertr\"aglich, d.h. es sind Situationen denkbar, in denen alle $B$ zwar
$A$ sind, aber
keines notwendigerweise. Diese Argumentati\-onsweise des Aristoteles
legitimiert auch das oben angef\"uhrte Gegenbeispiel, dass ein
Individuum unter einen Begriff mit Notwendigkeit f\"allt, unter einen
anderen aber nur zuf\"alligerweise. \"Ubrigens stimmt das Gegenbeispiel des
Aristoteles zu Barba\-ra XN, $A$ Bewegung, $B$ Sinnenwesen, $C$ Mensch, mit
dem oben genannten in einer Welt \"uberein, in der es genau ein Bewegtes,
ein Sinnenwesen und einen Menschen gibt.
\par\bigskip
\bigskip
\noindent
\unitlength0.5cm
\begin{picture}(10,4)
\put(0,3.3){\rm Celarent XN }
\put(0,2){$ \forall x (Bx \Longrightarrow  \neg Ax) $}
\put(0,1){$ \forall x(Cx \Longrightarrow NBx) $}
\put(13,3.3){\rm Ferio XN }
\put(13,2){$\forall x( Bx \Longrightarrow \neg Ax) $}
\put(13,1){$ \exists x(Cx \wedge NBx) $}
\end{picture}
\par\bigskip
\noindent
Ein gemeinsames Gegenbeispiel hierzu wird gegeben durch
$W_0=\{ x_0,y_0 \}$ und $W_1 = \{ x_1 \} $ mit den beiden
Individuen $x=(x_0,x_1)$ und $y=(y_0,x_1)$ und den Begriffen
$A_0=\{ y_0 \}$, $ B_0=C_0=\{ x_0 \}$, $A_1=B_1=C_1 = \{ x_1 \}$.
Damit ist $NAy,\, NBx ,\, NCx ,\, M(x =y) $. 
Die beiden Pr\"amissen sind erf\"ullt, 
die in Frage stehenden Konklusionen
$ \forall x \forall y(Cx \wedge Ay \Longrightarrow N(x \neq y))$
f\"ur Celarent bzw.
$\exists x[Cx \wedge \forall y (Ay \Longrightarrow N(x \neq y))]$
oder
$\exists x [NCx \wedge \neg Ax
\wedge  \forall y(NAy \Longrightarrow N(x \neq y))] $ f\"ur Ferio
dagegen nicht. Das Beispiel zeigt, dass noch nicht einmal bei der
starken partikul\"ar bejahenden Voraussetzung $\exists x (NCx \wedge NBx)$
ein Notwendigkeitsschluss stattfindet. Es gilt auch nicht $N(BoA)$ und
das Beispiel l\"asst sich auch als Gegenbeispiel zu Bocardo XN verwenden.
\newpage
\noindent
2. Figur
\par\bigskip
\noindent
Einige Notwendigkeitsschl\"usse der zweiten und der dritten Figur
lassen sich unter Ausnutzung der Symmetrie von $N(BeA)$ und $N(BiA)$
leicht auf Not\-wendigkeitsschl\"usse der ersten Figur zur\"uckf\"uhren.
Dieses bequeme Verfah\-ren, von dem auch Aristoteles immer wieder Gebrauch
macht, soll im Folgen\-den Verwendung finden und anhand einer
Ableitung f\"ur Camestres XNN ver\-deut\-licht werden.
\par\bigskip
\noindent
Camestres XNN
\par\bigskip
\noindent
\unitlength0.5cm
\begin{picture}(25,7)
\put(0,6){ $(1) \forall y(By \Longrightarrow Ay) $ }
\put(15,6){\rm  1. Pr\"amisse}
\put(0,4.8){ $(2)  \forall x \forall y (Cx \wedge Ay
\Longrightarrow N(x \neq y)) $}
\put(0,3.6){ $(3)  \forall y \forall x(Ay \wedge Cx
\Longrightarrow N(x \neq y)) $}
\put(0,2.4){ $(4)  \forall y \forall x(By \wedge Cx
\Longrightarrow N(x \neq y)) $}
\put(0,1.2){ $(5)  \forall x \forall y(Cx \wedge By
\Longrightarrow N(x \neq y))) $}
\put(15,4.8){\rm 2. Pr\"amisse}
\put(15,3.6){\rm Symmetrie von $N(CeA)$ auf (2)}
\put(15,2.4){\rm Celarent NXN auf (3) und (1)}
\put(15,1.2){\rm Symmetrie von $N(BeC)$ auf (4)}
\end{picture}
\par\bigskip
\noindent
\"Ahnlich l\"asst sich Cesare NXN und Festino NXN auf Celarent NXN bzw.
Ferio NXN zur\"uckf\"uhren und beweisen.
\par
In den \"ubrigen Situationen der zweiten Figur lassen sich keine
Notwen\-digkeitsschl\"usse ziehen.
\par\bigskip
\bigskip
\noindent
\unitlength0.5cm
\begin{picture}(25,4)
\put(0,3.3){\rm Camestres NX }
\put(0,2){$ \forall x (Bx \Longrightarrow NAx) $}
\put(0,1){$ \forall x(Cx \Longrightarrow  \neg Ax) $}
\put(13,3.3){\rm Cesare XN }
\put(13,2){$\forall x(Bx \Longrightarrow \neg Ax) $}
\put(13,1){$ \forall x(Cx \Longrightarrow NAx) $}
\end{picture}
\par\bigskip
\noindent
Diese Situationen kann man wieder auf Celarent XN zur\"uckf\"uhren oder
aber das Gegenbeispiel von dort unter Vertauschung der Rollen der
Mengen $A$,$B$ und $C$ hierher \"ubertragen. Auch Aristoteles erh\"alt die
erste Figur durch Um\-kehrung der verneinenden Pr\"amisse; ferner widerlegt er
f\"ur Camestres NX die An\-nahme, der Schlusssatz sei notwendig: ``Wenn
ferner der Schlusssatz not\-wendig ist, ergibt sich, dass $C$ einem $A$
notwendig nicht zukommt. Denn wenn $B$ notwendig keinem $C$ zukommt, wird
auch $C$ notwendig keinem $B$ zu\-kom\-men. Allein $B$ kommt notwendig einem $A$
zu, da ja auch $A$ notwendig jedem $B$ zukam. Und so muss denn $C$ einem
$A$ notwendig nicht zukommen.''\footnote{30b24ff}  
For\-malisiert bedeutet dies: Sei die in Frage stehende Konklusion
$N(CeB) \Longleftrightarrow
\forall x \forall y (Bx \wedge Cy \Longrightarrow N(x \neq y))$
als gegeben angenommen. Unter der Voraus\-setzung $B \neq \emptyset$
liefert die erste Pr\"amisse $ \exists x(NAx \wedge Bx)$. Diese
beiden Aus\-dr\"ucke ergeben dann kombiniert
$\exists x[NAx \wedge \forall y (Cy \Longrightarrow N(x \neq y))]$,
also (1) aus dem Implikationsdiagramm f\"ur $N(AoC)$ und damit insbesondere
$N(AoC)$. Aristoteles wendet dagegen ein: ``Aber es steht nichts
im Wege, $A$ so zu fassen, dass ihm seinem ganzen Umfange nach $C$ zukommen
kann.''\footnote{Erste Analytik 30b24ff}
(Obwohl es ihm tats\"achlich nicht zukommt), d.h.
die Pr\"amissen in Cames\-tres NX sind mit der Nebenbedingung
$\forall x(Ax \Longrightarrow \exists y(Cy \wedge M(x =y)))$
ver\-tr\"ag\-lich, die im Widerspruch zu dem erzielten Resultat steht.
\par\bigskip
\noindent
\unitlength0.5cm
\begin{picture}(10,4.5)
\put(0,3){\rm Festino XN }
\put(0,1.8){$ \forall x (Bx \Longrightarrow \neg Ax) $}
\put(0,0.6){$ \exists x(Cx \wedge  NAx) $}
\end{picture}
\par\bigskip
\noindent
Vertauscht man im Obersatz wieder $A$ und $B$, so ergibt sich Ferio XN,
wo kein Schluss auf Notwendigkeit stattfindet. Der letzte Schluss der
zweiten Figur, Baroco, l\"asst sich nicht auf eine Situation der ersten
Figur zur\"uckf\"uhren. Bei Baroco ergibt sich, wie sp\"ater bei Bocardo,
keinerlei Notwendigkeitsschluss.
\par\bigskip
\noindent
\unitlength0.5cm
\begin{picture}(15,4.5)
\put(0,3){\rm Baroco NX }
\put(0,1.8){$ \forall x (Bx \Longrightarrow NAx) $}
\put(0,0.6){$ \exists y(Cy \wedge  \neg Ay ) $}
\end{picture}
\par\bigskip
\noindent
Ein Gegenbeispiel erh\"alt man so. Es sei $W_0 = \{ x_0 , y_0 \} $ und
$W_1 = \{ x_1 \} $ mit den beiden Individualkonzepten
$x=(x_0,x_1),\, y=(y_0,x_1)$.
Es ist $x \neq y $, aber $M(x =y)$, da ja $y_1 =x_1 $.
Die Begriffe $A,B,C$ seien gegeben durch $A_0 =B_0 = \{ x_0 \} $,
$C_0 = \{ y_0 \} $ und durch $A_1 =B_1 = C_1= \{ x_1 \}$.
Damit ist $NAx$ und $NBx$
Damit sind die beiden Pr\"amissen erf\"ullt, dagegen
gilt noch nicht einmal
$\exists y(Cy \wedge \forall x (NBx \Longrightarrow N(x \neq y)))$,
also auch nicht $N(CoB)$.
\par\bigskip
\noindent
\unitlength0.5cm
\begin{picture}(25,4.5)
\put(0,3){\rm Baroco XN }
\put(0,1.8){$ \forall x (Bx \Longrightarrow Ax) $}
\put(0,0.6){$ \exists y[NCy  \wedge  \neg Ay \wedge 
\forall x (NAx \Longrightarrow N(x \neq y))] $}
\end{picture}
\par
\bigskip
\noindent
Hier folgt lediglich $ \exists y (NCy \wedge \neg By) $, das ist der
Ausdruck (5) im Implikations\-diagramm f\"ur $N(CoB)$.
Ein Schluss auf $N(CoB)$ kann sich dagegen
nicht ergeben, wie das folgende Gegenbeispiel lehrt:
$W_0=\{x_0 , y_0 \}$ und $W_1 = \{ x_1 \}$ mit den beiden Individuen
$x$ und $y$ mit $M(x =y)$ und den Begriffen
$A_0=B_0= \{ x \} \, , C_0 = \{ y_0 \}$ und
$A_1 = \emptyset \, , B_1 =C_1 =\{ x_1 \} $.
Es ist dann $ \neg NAx$ und $ NBx$, $NCy$.
Die Pr\"amissen sind erf\"ullt, der Ausdruck
$ (4) \exists y(Cy \wedge \neg By \wedge
\forall x(NBx \Longrightarrow N(x \neq y ))) $
aus dem Implikationsdiagramm aber
nicht, da ja $NBx$ gilt, aber eben nicht $N(x \neq y)$.
Somit gilt auch nicht $N(CoB)$.
\par
\bigskip
\bigskip
\noindent
3. Figur
\par\bigskip
\noindent
Der erste Schluss der dritten Figur, Darapti, ist der einzige der
vierzehn Schl\"usse, bei dem unabh\"angig davon, welche Pr\"amisse notwendig
ist, die Konklusion notwendig ist. Allerdings setzt dies voraus,
da\ss\ $C$ nicht leer ist, eine Annahme, von der Aristoteles bereits
im assertorischen Fall Darapti XXX stillschweigend Gebrauch macht.
Unter dieser Voraussetzung l\"asst sich jeweils die assertorische
Allaussage auf eine assertorische Existenzaussage reduzieren, so dass sich
der Notwendigkeitsschluss mittels Darii NXN ergibt.
\par\bigskip
\bigskip
\noindent
Darapti NXN ($C \neq \emptyset$)
\par\bigskip
\noindent
\unitlength0.5cm
\begin{picture}(25,7)
\put(0,6){(1) $ \forall x(Cx \Longrightarrow NAx) $ }
\put(17,6){\rm  1. Pr\"amisse}
\put(0,4.8){(2) $  \forall x (Cx \Longrightarrow Bx) $}
\put(0,3.6){(3) $  \exists x (Bx \wedge Cx) $}
\put(0,2.4){(4)  $\exists x(Bx \wedge NAx)$}
\put(0,0.9){\rm also $N(BiA)$}
\put(17,4.8){\rm 2. Pr\"amisse}
\put(17,3.6){\rm Existenzreduzierung auf (2)}
\put(17,2.4){\rm Darii NXN auf (1) und (3)}
\end{picture}
\par\bigskip
\noindent
Ebenso ergibt sich bei Darapti XNN ($C \neq \emptyset$) $N(BiA)$. \"Ahnlich
verh\"alt es sich mit Felapton NXN; auch hier reduziert man unter der
Voraussetzung, dass $C$ nicht leer ist, die allgemein bejahende zweite
Pr\"amisse zu einer Existenzaussage und schlie\ss t mit Ferio NXN.
\par\bigskip
\noindent
Felapton NXN ($C \neq \emptyset$)
\par\bigskip
\noindent
\unitlength0.5cm
\begin{picture}(10,7)
\put(0,6){(1) $ \forall x \forall y( Cx \wedge Ay
\Longrightarrow N(x \neq y)) $ }
\put(17,6){\rm  1. Pr\"amisse}
\put(0,4.8){(2) $  \forall x (Cx \Longrightarrow Bx) $}
\put(0,3.6){(3) $  \exists x (Bx \wedge Cx) $}
\put(0,2.4){(4)  $\exists x(Bx \wedge
\forall y (Ay \Longrightarrow N(x \neq y))$}
\put(0,0.9){\rm also $N(BoA)$}
\put(17,4.8){\rm 2. Pr\"amisse}
\put(17,3.6){\rm Existenzreduzierung auf (2)}
\put(17,2.4){\rm Ferio NXN auf (1) und (3)}
\end{picture}
\par\bigskip
\noindent
Die beiden folgenden Schl\"usse lassen sich leicht auf Darii NXN
zur\"uckf\"uhren, indem man die Symmetrie des partikul\"ar bejahenden
ausn\"utzt; man kann aber auch direkt schlie\ss en.
\par\bigskip
\bigskip
\noindent
\unitlength0.5cm
\begin{picture}(10,6)
\put(0,5.3){\rm Datisi NXN }
\put(0,4){$ \forall x (Cx \Longrightarrow NAx) $}
\put(0,3){$ \exists x(Cx \wedge Bx) $}
\put(0,2.6){\---------------------------}
\put(0,1.7){$ \exists x(NAx \wedge Bx)$}
\put(13,5.3){\rm Disamis XNN}
\put(13,4){$\exists x(Cx \wedge Ax) $}
\put(13,3){$ \forall x(Cx \Longrightarrow NBx) $}
\put(13,2.6){\-----------------------}
\put(13,1.7){$ \exists x(NBx \wedge Ax) $}
\end{picture}
\par
\noindent
Schlie\ss lich ist Ferison NXN \"aquivalent zu Ferio NXN.
\par
Keine Notwendigkeitsschl\"usse finden hingegen in den \"ubrigen Situatio\-nen
der dritten Figur statt.
\par\bigskip
\bigskip
\noindent
\unitlength0.5cm
\begin{picture}(10,5)
\put(0,4.3){\rm Felapton XN }
\put(0,3){$ \forall x (Cx \Longrightarrow \neg Ax) $}
\put(0,2){$ \forall x(Cx \Longrightarrow NBx) $}
\put(13,4.3){\rm Ferison XN }
\put(13,3){$\forall x(Cx \Longrightarrow \neg Ax) $}
\put(13,2){$ \exists x(Cx \wedge NBx) $}
\end{picture}
\par\smallskip
\noindent
Hier folgt (bei $C \neq \emptyset $) lediglich
$ \exists x (NBx \wedge \neg Ax)$, also der
Ausdruck (5) im Implikations\-diagramm, der unterhalb von $N(BoA)$
angesiedelt ist. Das Ge\-genbeispiel zu Celarent XN / Ferio XN
l\"asst sich hierher
\"ubertragen. Ebenso kann man das Gegenbeispiel zu Darii XN unter
Vertau\-schung der Rollen von $B$ und $C$ verwenden, um zu zeigen, dass
sich bei Datisi XN und Disamis NX keine notwendige Konklusion ergibt.
\par
Es bleibt noch zu zeigen, dass sich bei Bocardo kein Notwendigkeits\-schluss
ziehen l\"asst, unabh\"angig davon, welche der beiden Pr\"a\-mis\-sen
apodiktisch ist.
\par\bigskip
\noindent
\unitlength0.5cm
\begin{picture}(10,4.5)
\put(0,3){\rm Bocardo XN }
\put(0,1.8){$ \exists x (Cx  \wedge \neg Ax) $}
\put(0,0.6){$ \forall x (Cx \Longrightarrow NBx) $}
\end{picture}
\par\bigskip
\noindent
Hier folgt wieder nur $\exists x (NBx \wedge \neg Ax)$. Das
Gegenbeispiel zu Ferio XN dient hier als solches.
\par\bigskip
\noindent
\unitlength0.5cm
\begin{picture}(10,4.5)
\put(0,3){\rm Bocardo NX }
\put(0,1.8){$ \exists x[NCx \wedge \neg Ax \wedge
\forall y (NAy \Longrightarrow N(x \neq y))]$}
\put(0,0.6){$ \forall x (Cx \Longrightarrow Bx) $}
\end{picture}
\par\bigskip
\noindent
Hier folgt lediglich $\exists x[Bx \wedge \neg Ax \wedge
\forall y (NAy \Longrightarrow N(x \neq y))]$, d.i. der
Ausdruck (4) im Implikationsdiagramm f\"ur $N(BoA)$. Weder l\"asst sich
vorne $NBx$ ab\-leiten noch l\"asst sich hinten $NAy$ zu
$Ay$ verallgemeinern. Ein Gegenbei\-spiel ist 
gegeben durch $W_0=\{ x_0,y_0 \} $ und $ W_1 = \{ x_1 \} $, den beiden
Individuen $x,y$ und den Begriffen
$A_0 =\{ y_0 \}, \, B_0=C_0=\{ x_0 \}$ und
$C_1 =\{ x_1 \} , \, A_1=B_1= \emptyset \, .$
\par\bigskip
\noindent
Es steht noch aus, den Nachweis der Eigenschaft (4) einer vollst\"andigen
Formalisierung zu erbringen, also zu zeigen, dass bei s\"amtlichen
Schl\"ussen die Konklusion notwendig ist, wenn dies f\"ur beide
Pr\"amissen gilt (Kap.8). Wie aus den Kapiteln 9-11 hervorgeht und
wie oben gezeigt worden ist, sind Baroco und Bocardo die einzigen
Schl\"usse, bei denen, gleichg\"ultig ob der Ober- oder der Untersatz
notwendig ist, kein Schluss auf Notwendigkeit stattfindet. F\"ur die
\"ubrigen zw\"olf Situationen findet zumindest in einem Fall, \--- wenn
der Obersatz oder der Untersatz
notwendig ist \--- ein Not\-wendigkeits\-schluss
statt. Sind in diesen zw\"olf Situationen nun beide Pr\"amissen
notwen\-dig, so ergibt sich die notwendige Konklusion einfach dadurch,
dass man eine geeignete apodiktische Pr\"amisse auf ihr assertorisches
Analogon zur\"uck\-f\"uhrt, was aufgrund des $T-$Axioms bzw. der Eigenschaft
(1) m\"oglich ist. So l\"asst sich beispielsweise Barbara NN zu Barbara NX
abschw\"achen und von da aus die notwendige Konklusion erschlie\ss en.
\par
F\"ur Baroco NNN und Bocardo NNN kann eine solche Argumentation nat\"urlich
nicht funktionieren, und es verwundert nicht, dass diese auch bei
Aristoteles eine Sonderrolle spielen: ``(Bei diesen) wird der Beweis
nicht auf die gleiche Art gef\"uhrt werden
k\"onnen, (...)''\footnote{Erste Analytik, 30a6ff}. 
Ein Beweis ergibt sich nun wie folgt,
wobei es sogar gen\"ugt, f\"ur $N(CoA)$ den nach dem
Implikationsdia\-gramm schw\"acheren Ausdruck (4) anzusetzen, um die
gew\"unschte Konklusion $N(CoB)$ bzw. $N(BoA)$ zu erhalten.
\par\bigskip
\noindent
\unitlength0.5cm
\begin{picture}(10,4.5)
\put(0,3.8){\rm Baroco NNN }
\put(0,2.5){$ \forall y (By \Longrightarrow NAx) $}
\put(0,1.5){$ \exists x[Cx \wedge \neg Ax  \wedge
\forall y (NAy \Longrightarrow N(x \neq y))] $}
\put(0,1.1){\------------------------------------------------------}
\put(0,0.2){$ \exists x[Cx \wedge
\forall y(By \Longrightarrow N(x \neq y))]$}
\end{picture}
\par\bigskip
\noindent
F\"ur das $x$ aus dem Untersatz folgt ja f\"ur ein beliebiges $y \in B$ mit
dem Obersatz $NAy$ und daraus mit dem Untersatz $N(x \neq y)$. So
erh\"alt man also (2) aus dem Implikationsdiagramm und damit auch $N(CoB)$.
\par\bigskip
\noindent
\unitlength0.5cm
\begin{picture}(10,4.5)
\put(0,3.8){\rm Bocardo NNN }
\put(0,2.5){$ \exists x [Cx \wedge \neg Ax \wedge
\forall y(NAy \Longrightarrow N(x \neq y))] $}
\put(0,1.5){$ \forall x(Cx \Longrightarrow NBx) $}
\put(0,1.1){\------------------------------------------------------}
\put(0,0.2){$ \exists x[NBx \wedge \neg Ax 
\wedge \forall y(NAy \Longrightarrow N(x \neq y))]$}
\end{picture}
\par\bigskip
\noindent
Es ergibt sich also der Ausdruck (3) aus dem Implikationsdiagramm
und damit auch $N(BoA)$.
\par
\bigskip
\bigskip
\bigskip
\noindent
IV. ZU DEN KONTIGENTEN SCHL\"USSEN DES ARISTOTELES
\par\bigskip
\noindent
Zum Schluss m\"ochte ich andeuten, wie die vorgestellte modallogische
Inter\-pretation auf die kontigenten Schl\"usse des Aristoteles ausgedehnt
werden kann. Dabei ist nicht an eine vollst\"andige Behandlung
s\"amtlicher Schl\"usse der Kapitel 14-22 gedacht, sondern lediglich an
die Untersuchung einiger charakteristischer Situationen. Ich beschr\"anke
mich auf die allgemeinen Kontingenzbeziehungen und auf die Schl\"usse
Barbara, Celarent, Cesare und Camestres.
\par
Vom kontingenten Zukommen fordert Aristoteles, dass ``(...) es
ist kon\-tin\-gent zu sein,
mit dem Satz: es ist kontingent nicht zu sein, vertauscht
werden kann (...)''.\footnote{32a32ff}
Hiermit ergibt sich dann die immer wieder verwendete
\"Aquivalenz von ``$A$ kann jedem $B$ zukommen'' und ``$A$ kann jedem $B$ nicht
zukommen''. Setzt man f\"ur die allgemeinen Kontingenzbeziehungen
Aus\-dr\"ucke der Form
\par\bigskip
\noindent
\unitlength0.5cm
\begin{picture}(10,3)
\put(0,2.4){$ K(BaA) : \Longleftrightarrow
\forall x ( Bx \Longrightarrow KAx ) $ }
\put(0,1){$ K(BeA) :\Longleftrightarrow
\forall x ( Bx \Longrightarrow K( \neg Ax)) \,\, , $}
\end{picture}
\par
\noindent
so ergibt sich die von Aristoteles beschriebene \"Aquivalenz
$K(BaA) \Longleftrightarrow K(AeB)$ am einfachsten mittels der \"Aquivalenz
$KAx \Longleftrightarrow K(\neg Ax)$. Dieses elementare kontingente
Zukommen muss dann durch einen pr\"adikatenlo\-gischen Ausdruck
wiedergegeben werden, der sowohl einen bejahenden als auch einen
verneinenden Teilaspekt besitzt.
Welche Mikrostruktur f\"ur $KAx$ geeignet
ist wird weiter unten untersucht, wenn die kontingent\--assertorischen
und die kontingent\--apodiktischen Schl\"usse behandelt werden.
F\"ur die reinen Kontingenzschl\"usse ist diese Struktur ohne Belang,
stattdessen dr\"angt sich eine andere Be\-mer\-kung des Aristoteles in den
Mittelpunkt der \"Uberlegung. Im 13. Kapitel schreibt er: ``Denn wenn
man sagt, wovon $B$, davon kann $A$ ausgesagt werden, so hat das eine von
diesen beiden Be\-deu\-tung\-en: ent\-we\-der wo\-von $B$ aus\-ge\-sagt wird,
oder wo\-von es aus\-ge\-sagt wer\-den
kann;''\footnote{32b25ff} Aristoteles unterscheidet
also zwischen $\forall x(Bx \Longrightarrow KAx)$ und
$\forall x(KBx \Longrightarrow KAx )$,
ohne dass v\"ollig klar wird, wann er mit welcher Kontingenzbeziehung
zu arbeiten beabsichtigt. Seine weiteren Ausf\"uhrungen legen aber
den Schluss nahe, dass die zwei\-te Be\-zie\-hung zumin\-dest im sich
un\-mit\-tel\-bar an\-schlie\ss \-en\-den Ka\-pi\-tel 14 ver\-wendet wird.
Es ergeben sich dann leicht die rein kontin\-genten
Schl\"usse der ersten Figur.
\par\bigskip
\noindent
\unitlength0.5cm
\begin{picture}(10,4.5)
\put(0,3.8){\rm Barbara KKK}
\put(0,2.5){$ \forall x (KBx \Longrightarrow KAx ) $}
\put(0,1.5){$ \forall x(KCx \Longrightarrow KBx) $}
\put(0,1.1){\----------------------------}
\put(0,0.2){$ \forall x(KCx \Longrightarrow KAx )$}
\end{picture}
\par\bigskip
\noindent
Es verwundert nicht, dass Aristoteles diesen Schluss als ``vollkommen''
be\-zeichnet. Als unvollkommen gilt ihm dagegen ein Schluss wie
\par\bigskip
\noindent
\unitlength0.5cm
\begin{picture}(10,3.2)
\put(0,2.5){$ \forall x (KBx \Longrightarrow KAx) $}
\put(0,1.5){$ \forall x(KCx \Longrightarrow K(\neg Bx)) $}
\put(0,1.1){\---------------------------------}
\put(0,0.2){$ \forall x(KCx \Longrightarrow KAx)$}
\end{picture}
\par\bigskip
\noindent
da hierzu die \"Aquivalenz $K(\neg Bx) \Longleftrightarrow KBx$
verwendet, also eine ``Um\-keh\-rung''  vorgenommen wird. 
\par\bigskip
\noindent
Wenden wir uns nun den interessanteren und schwierigeren
kon\-tin\-gent\--as\-sertorischen und kon\-tingent\--apodiktischen Schl\"ussen zu.
Es empfiehlt sich, von den beidseitigen Kontingenzbeziehungen Abschied
zu nehmen und im Folgenden stets die einseitige Kontingenzbeziehung,
also $\forall x(Bx \Longrightarrow KAx)$ anzusetzen. Des weiteren wird es
unumg\"anglich sein, auf die Mikrostruktur von
$KAx$ einzugehen, was aber f\"ur einen kurzen Moment
noch zur\"uckgestellt werden kann, da sie bei einigen Schl\"ussen keine
Rolle spielt, wie beispielsweise in der folgenden, von Aristoteles als
vollkom\-men bezeichneten Situation.
\par\bigskip
\noindent
\unitlength0.5cm
\begin{picture}(10,4.5)
\put(0,3.8){\rm Barbara KXK (ebenso Celarent KXK)}
\put(0,2.5){$ \forall x (Bx \Longrightarrow KAx) $}
\put(0,1.5){$ \forall x(Cx \Longrightarrow Bx) $}
\put(0,1.1){\----------------------------}
\put(0,0.2){$ \forall x(Cx \Longrightarrow KAx)$}
\end{picture}
\par\bigskip
\noindent
Damit hat nun der Kontingenzoperator $K$ seine Aufgabe als Platzhalter
erf\"ullt und es gilt, f\"ur $KAx $ bzw. $K(\neg Ax)$ einen
Ausdruck in der bisher verwendeten Sprache der Modallogik anzugeben, der
einerseits die bereits erw\"ahnte und verwendete \"Aquivalenz
$KAx \Longleftrightarrow  K(\neg Ax)$ einleuchtend macht
und andererseits gestattet, die weiteren von
Aristoteles gezogenen Schl\"usse nachzuzeichnen. Der erste
Punkt legt nahe, dass die Feinstruktur von $KAx$ sowohl einen
m\"oglicherweise\--bejahenden als auch einen
m\"oglicherweise\--ver\-neinenden
Teil aufweisen muss, also von der Gestalt
$\bar{M}Ax \wedge \tilde{M}( \neg Ax)$ ist. F\"ur
diese elementaren M\"oglichkeitsbeziehungen sind wiederum verschiedene
Formulierungen denkbar, die innerhalb unseres modal\--pr\"adikatenlogischen
Rahmens bzw. in unserem Modell interpretierbar sind.
\par\bigskip
\noindent
\unitlength0.5cm
\begin{picture}(10,6.5)
\put(0,5.7){\rm Kontingentes Zukommen $\bar{M}Ax$}
\put(0,4){$ (1) \,\, MAx $}
\put(0,2.6){$ (2)\,\, \exists y(Ay \wedge M(x =y)) $}
\put(0,1.2){$ (3)\,\, \exists y(MAy \wedge M(x=y))$}
\put(13,5.7){\rm Kontingentes Nichtzukommen $\tilde{M}(\neg Ax) $}
\put(13,4){$(1)\,\, M(\neg Ax ) $}
\put(13,2.6){$(2)\,\,  \exists y (\neg Ay \wedge M(x=y)) $}
\put(13,1.2){$(3)\,\,  \exists y(M(\neg Ay) \wedge M(x=y)) $}
\end{picture}
\par
\noindent
Dabei folgen all diese kontingenten (Nicht\--)Zugeh\"origkeiten aus den
ent\-sprechenden assertorischen Beziehungen.
$(1)$ und $(2)$ sind implikativ nicht verbunden, wie man sich
in einem Modell leicht klar macht, und der Ausdruck $(3)$ folgt
sowohl aus $(1)$ als auch aus $(2)$.
\par
Welcher dieser Ausdr\"ucke ist nun am ehesten geeignet, die Pl\"atze in
$KAx \Longleftrightarrow K(\neg Ax) \Longleftrightarrow
\bar{M}Ax \wedge \tilde{M} (\neg Ax) $
einzu\-nehmen? Auf\-grund den
Erfahr\-ungen bei der Behand\-lung der apodik\-tischen Schl\"usse ist kaum
zu erwarten, dass f\"ur den positiven und den nega\-tiven Teil analoge
Formulie\-rungen einzusetzen sind. Daher kann man sich zun\"achst nur auf
den beja\-hen\-den Teil konzentrieren und unter\-suchen, wie dieser gestaltet
sein muss, um mit den aristotelischen Er\-geb\-nissen \"Ubereinstimmung
zu erzielen. Es sei daran erinnert, dass die fol\-genden
Schl\"usse keine Konklusion in dem Sinne haben, dass etwas
einem zugleich m\"oglicherweise und m\"oglicherweise nicht zukommt,
sondern nur in dem Sinne, dass etwas einem m\"oglich\-erweise
zukommt ohne ihm unbedingt auch m\"oglich\-erweise nicht zukom\-men
zu k\"on\-nen oder umge\-kehrt.
\par\bigskip
\noindent
\unitlength0.5cm
\begin{picture}(10,4)
\put(0,3.3){\rm Barbara XK}
\put(0,2){$ \forall x (Bx \Longrightarrow Ax) $}
\put(0,1){$ \forall x(Cx \Longrightarrow
\bar{M}Bx \wedge \tilde{M} \neg Bx) $}
\end{picture}
\par
\noindent
Aristoteles schlie\ss t hieraus, dass alle $C$ m\"oglicherweise $A$ sind.
Seine Aus\-f\"uh\-rungen im Vorfeld dieses Schlusses lassen erkennen,
dass er hier eine Mo\-no\-tonieregel gelten l\"asst.
Die bejahende kontingente Aussage ist also nicht nur im Subjekt, sondern
auch im Pr\"adikat extensional durchl\"assig. Eine Monotonieregel, die
erlaubt, aus $B \subseteq A $ und $MBx$ auf $MAx$ zu
schlie\ss en, gilt in unseren modallogischen Modellen nicht und w\"urde
die gesamte Formali\-sierung der Notwendigkeitssyllogistik
zunichte machen. Von daher bleibt nur \"ubrig zu setzen
$\bar{M}Ax \Longleftrightarrow (2)\, \exists y(Ay \wedge M(x=y))$,
und dieser Ausdruck ist aufgrund der au\ss ermodalen Stellung
des Pr\"adikats automatisch, also ohne eine zus\"atzliche Regel monoton.
Damit erhalten wir\footnote{Hier und im Folgenden werden in den
Kontingenzbeziehungen nur der Teil angef\"uhrt, der f\"ur den Schluss
von Bedeutung ist}
\par\bigskip
\noindent
\unitlength0.5cm
\begin{picture}(10,4.5)
\put(0,3.8){\rm Barbara XKM}
\put(0,2.5){$ \forall y (By \Longrightarrow Ay) $}
\put(0,1.5){$ \forall x(Cx \Longrightarrow
\exists y(By \wedge M(x=y)))$}
\put(0,1.1){\---------------------------------------------}
\put(0,0.2){$ \forall x(Cx \Longrightarrow
\exists y(Ay \wedge M(x=y)))$}
\end{picture}
\par\bigskip
\noindent
Wir fixieren
$\bar{M}Ax \Longleftrightarrow \exists y( Ay \wedge M(x=y))$,
und dieser Ausdruck wird auch in den \"ubrigen Situationen die
gew\"unschten Resultate erbringen. Insbesondere erlaubt dieser Ausdruck
in Kapitel 16, wo die eine Pr\"amisse not\-wendig und die andere kontingent ist,
in den von Aristoteles be\-schrie\-benen Situa\-tionen den Schluss auf
assertorisches Nichtsein, wie anhand eines Be\-weises f\"ur Cela\-rent NKX
gezeigt werden soll.
\par\bigskip
\noindent
\unitlength0.5cm
\begin{picture}(10,4.5)
\put(0,3.8){\rm Celarent NKX}
\put(0,2.5){$ \forall y \forall x(By \wedge Ax \Longrightarrow N(x \neq y))$}
\put(0,1.5){$ \forall x(Cx \Longrightarrow \exists y(By \wedge M(x=y)))$}
\put(0,1.1){\---------------------------------------------}
\put(0,0.2){$ \forall x(Cx \Longrightarrow \neg Ax )$}
\end{picture}
\par\bigskip
\noindent
Die Annahme $ \exists x(Cx \wedge Ax) $ f\"uhrt mit der zweiten
Pr\"amisse zu $\exists x \exists y(Ax \wedge By \wedge M(x=y))$, was gerade
durch die erste Pr\"amisse ausgeschlossen wird.
\par
Sind dagegen beide Pr\"amissen bejahend, wie bei Barbara KNK, so ergibt sich
kein assertorischer Schlusssatz. Reduziert man die notwendige
Aussage auf ihr assertorische Analogon, so erh\"alt man die schon behandelten 
kontin\-gent\--assertorischen Situationen und damit Schl\"usse auf Kontingenz.
\par\bigskip
\noindent
Wenden wir uns nun dem negativen Teil der M\"oglichkeitsaussage zu.
Aus\-gehend von der Formalisierung des positiven Teils gibt Celarent XKM
ersten Auf\-schluss dar\"uber, wie dieser gestaltet sein muss. 
Bei Celarent XKM erh\"alt Aristoteles, dass m\"oglicherweise
$A$ keinem $C$ zukommt.
\par\bigskip
\noindent
\unitlength0.5cm
\begin{picture}(10,4.5)
\put(0,3.8){\rm Celarent XKM}
\put(0,2.5){$ \forall y (By \Longrightarrow \neg Ay) $}
\put(0,1.5){$ \forall x(Cx \Longrightarrow
\exists y(By \wedge M(x=y)))$}
\put(0,1.1){\----------------------------------------------}
\put(0,0.2){$ \forall x(Cx \Longrightarrow
\exists y(\neg Ay \wedge M(x=y)))$}
\end{picture}
\par\bigskip
\noindent
Als Konklusion erhalten wir, dass jedes Element aus $C$ gem\"a\ss\ 
dem Aus\-druck $(2)$ m\"oglicherweise nicht zu $A$ geh\"ort. Der
gesuchte Ausdruck f\"ur das m\"oglicherweise Nichtzukommen
$\tilde{M}(\neg Ax)$ muss also aus
$\exists y (\neg Ay \wedge M(x=y)) $
folgen. Die Gleichsetzung $\tilde{M}(\neg Ax) \Longleftrightarrow (2)$
verbietet sich aber, da nach Aristoteles bei Camestres XK kein Schluss
stattfindet und $\tilde{M}(\neg Ax)$ nicht wie der positive Teil monoton
sein darf.
\par\bigskip
\noindent
\unitlength0.5cm
\begin{picture}(10,4.5)
\put(0,3.8){\rm Camestres XK}
\put(0,2.5){$ \forall y (By \Longrightarrow Ay) $}
\put(0,1.5){$ \forall x(Cx \Longrightarrow
\exists y(\neg Ay \wedge M(x=y)))$}
\put(0,1.1){\----------------------------------------------}
\put(0,0.2){$ \forall x(Cx \Longrightarrow
\exists y(\neg By \wedge M(x=y)))$}
\end{picture}
\par\bigskip
\noindent
Auch der nicht-monotone Ausdruck $M(\neg Ax)$ verbietet sich.
Einer\-seits folgt er nicht aus $(2)$ und andererseits liefert er
bei Cesare KN im Gegen\-satz zur ausf\"uhrlichen 
Argu\-men\-tation des Aristo\-teles im 19.Kapitel einen asserto\-rischen Schluss.
Diese Argumentation besagt zuge\-spitzt, dass das notwendige Zukommen dem
m\"oglich\-erweise Nicht\-zukommen (wie es als negativer Teil\-aspekt des
m\"oglich\-erweise Zukommens verwendet wird) nicht widerspricht.
\par\bigskip
\noindent
\unitlength0.5cm
\begin{picture}(10,4.5)
\put(0,3.8){\rm Cesare KN}
\put(0,2.5){$ \forall x (Bx \Longrightarrow M(\neg Ax)) $}
\put(0,1.5){$ \forall x(Cx \Longrightarrow NAx)$}
\put(0,1.1){\-----------------------------}
\put(0,0.2){$ C \cap B = \emptyset $}
\end{picture}
\par\bigskip
\noindent
Von den drei oben angef\"uhrten Kandidaten f\"ur $\tilde{M}(\neg Ax)$
bleibt nur noch der Ausdruck $(3)\,  \exists y(M(\neg Ay) \wedge M(x=y)) $
\"ubrig. Dieser folgt aus (2) und erlaubt damit Celarent XKM. Ferner
ist er nicht monoton, ein
Schluss bei Cames\-tres XK ist also nicht m\"oglich. Ebensowenig erlaubt
er irgendeinen Schluss bei Cesare, insbesondere nicht den assertorischen
Schluss auf Disjunktheit: etwas kann notwendigerweise ein $A$ sein, aber
trotzdem m\"oglicherweise gleich einem etwas, das kein $A$ sein muss.
\par
Ein Problem ergibt sich allerdings bei Camestres NK. Hierzu \"au\ss ert
sich Aristoteles nicht explizit im Gegensatz zu Cesare KN, wo er jeden
Schluss bestreitet. Wenn die negative Pr\"amisse aber kontingent ist, so
ist wegen der Nichtumkehrbarkeit der negativen Kontingenzaussage nicht wie
bei den assertorischen oder apodiktischen Situationen mit Cesare schon
Camestres mit abgehandelt. Bei Camestres NK erhalten wir, wenn wir $(3)$
ansetzen, einen Schluss auf das kontingente Nichtsein. Die
Notwendigkeitsaussage $\forall y(By \Longrightarrow NAy)$
l\"asst sich kontraponieren zu
$\forall y( M (\neg Ay) \Longrightarrow \neg By)$,
was wir gleich als erste Pr\"amisse anf\"uhren.
\par\bigskip
\noindent
\unitlength0.5cm
\begin{picture}(10,4.5)
\put(0,3.8){\rm Camestres NKM}
\put(0,2.5){$ \forall y(M(\neg A y) \Longrightarrow \neg By)$}
\put(0,1.5){$ \forall x (Cx \Longrightarrow
\exists y (M(\neg Ay) \wedge M(x=y)))$}
\put(0,1.1){\-----------------------------------------------------}
\put(0,0.2){$ \forall x (Cx \Longrightarrow
\exists y(\neg By \wedge M(x=y)))$}
\end{picture}
\par\bigskip
\noindent
Wir erhalten hier, dass alle $C$ im Sinne von $(2)$ und damit auch im Sinne
von $(3)$ m\"oglicherweise nicht $B$ sind. Zwar zieht Aristoteles diesen
Schluss nicht, er behauptet aber auch nirgendwo dessen Ung\"ultigkeit.
Ungeachtet dieser kritischen Stelle scheint mir der Ausdruck
$$ KAx \Longleftrightarrow \, \exists y(Ay  \wedge M(x=y)) \wedge
\exists z(M(\neg Az) \wedge M(x=z)) $$
am besten geeignet, die
vorgestellte Formalisierung auf die kontingenten Schl\"usse auszudehnen.
\par\bigskip
\noindent

\end{document}